\documentclass{amsart}

\usepackage{amsthm,amsmath}

\usepackage[round]{natbib}

 \theoremstyle{remark} 
 \newtheorem*{prop}{Proposition}

\begin{document}

\title{Failure of Calibration is Typical}

\author{Gordon Belot} 
\address{Department of Philosophy, University of Michigan} 
\email{belot@umich.edu}
\thanks{Thanks to an anonymous referee, whose comments on an earlier version led to improvements. \\ Forthcoming in \emph{Statistics and Probability Letters.} \\ Version of May 7, 2013}

\begin{abstract} 
\citet{Schervish:1985uo} showed that every forecasting system is noncalibrated for uncountably many data sequences that it might see. This result is strengthened here: from a topological point of view, failure of calibration is typical and calibration rare. Meanwhile, Bayesian forecasters are certain that they are calibrated---this invites worries about the connection between Bayesianism and rationality.
\end{abstract}

\maketitle

\section{Introduction}

 Consider a context in which Nature presents us with an infinite binary sequence, one bit at a time. A \emph{probabilistic forecasting system} is a function that takes as input finite binary strings (records of data seen to date) and gives as output numbers in the unit interval (the forecast probability that the next bit seen will be a one). Following custom, we specialize to the case of weather forecasting systems: the data sequences encode, for each of some sequence of days, whether or not it snows in a given locale on that day (1=snow; 0=no snow); each day, the forecasting system is asked to output a probability for snow on that day, given the pattern of snowy and snowless days in the extant data set. 

\emph{Notation.} We denote data sequences by $\boldsymbol{a}$=($a_1,$ $a_2,$ $a_3,$ \ldots) and the probability output by a forecaster on day $k,$ based on knowledge of $a_1,$ \ldots $a_{k-1},$ by $\pi_k.$  We call $d_k= a_k - \pi_k$ the forecaster's \emph{discrepancy} on day $k.$ 

Motivated rivals of a forecaster will keep an eye on long-term behaviour of the forecaster's mean discrepancy---if this quantity does not tend to zero when specialized to days of a given type (all days, even-numbered days, days on which the forecast probability exceeds .9, \ldots) then the forecaster can be discredited as being systematically excessively optimistic or pessimistic concerning the chance of snow on such days (or, possibly, as oscillating permanently between excessive pessimism and optimism).    

This motivates the following definitions \citep[\S 5]{Dawid:1985qd}. A \emph{selection rule} $\delta$ is a function that takes as input a forecasting system $F$ and a data sequence $\boldsymbol{a}$ and gives as output a subsequence of the sequence of days---with whether or not day $k$ appears in the subsequence $\delta(F, \boldsymbol{a})$ depending only on $a_1,$ \ldots $a_{k-1}$ and on the forecast probability $\pi_k$ for day $k.$ A forecaster $F$ is \emph{calibrated} for data sequence $\boldsymbol{a}$ with respect to selection rule $\delta$ if the mean discrepancy tends to zero when one specializes to the days picked out by $\delta(F,\boldsymbol{a})$ (or if there are only finitely many such days).

Each set of selection rules determines a notion of calibration.  Here we will just need to mention two, the first very stringent indeed, the other so weak that it will be implied by just about any interesting notion of calibration. (i) A forecasting system is \emph{computably calibrated} for data sequence $\boldsymbol{a}$ if it is calibrated for $\boldsymbol{a}$ relative to each computable selection rule. (ii)  A forecasting system is high-low calibrated for data sequence $\boldsymbol{a}$ if it is calibrated relative to the rule that selects those days on which $\pi_k \geq .5$ and the rule that selects those days on which $\pi_k < .5.$  

It is natural to ask whether there exist forecasting systems that are calibrated (in one or another sense) for each possible data sequence. Answer: there is no forecasting system that is even high-low calibrated for every possible data sequence \citep{Oakes:1985dn, Dawid:1985mz}. Fix a forecasting system and think of Nature as an adversary aiming to spoil the forecaster's high-low calibration. Nature can achieve this by basing the weather each day on the forecast probability of snow, making the day snowy if and only if the forecast probability is less than .5.\footnote{What if we allow the forecaster to play a mixed strategy, randomly selecting on each day a deterministic forecasting scheme whose forecast to follow? In itself, this will make no difference---all that Nature cares about is the forecast probability output each day \citep[\S 3]{Foster:1998rt}. But there are some interesting variants of our game that forecasters playing mixed strategies \emph{can} win: see, e.g., \citet{Foster:1998rt}, \citet{Vovk:2005gf}, and \citet{Vyugin:2009rz}.} 

 \citet{Schervish:1985uo} showed more: every forecasting system fails to be high-low calibrated relative to an uncountable infinity of possible data sequences;  indeed there are uncountably many data sequences with the feature that each computable forecasting method fails to be high-low calibrated relative to each one of them. He concluded that ``noncalibrability is as much the normal state as calibrability.'' 
 
 One way to build a forecasting system is to choose a prior over the space of binary sequences, and to use a Bayesian agent who conditionalizes that prior on finite data sets to output the desired probabilities. \citet{Dawid:1982rc, Dawid:1985qd} proved that every such Bayesian forecasting system is supremely confident of its own calibration: every prior on the space of binary sequences considers the set of sequences relative to which it is \emph{not} computably calibrated to be a set of measure zero.  \citet[Note 7.3]{Dawid:1985qd} saw in this a reply to Schervish: there remains an intuitive sense in which those sequences for which no computable forecasting system is high-low calibrated must be sparse---since they are, collectively, assigned measure zero by any computable Bayesian forecasting system.  

Now, it is a familiar point that there are uncountable subsets of the real line that are nonetheless small or sparse in an interesting sense: the ternary Cantor set, for instance, is uncountable but is a set of Lebesgue measure zero and is also meagre  (negligible from the topologists' point of view).  One purpose of the present note is to extend Schervish's result (and to undercut Dawid's response) by showing that if we consider topology rather than just cardinality, then the parity between calibration and noncalibration is broken: each forecasting system is fails to be high-low calibrated relative to typical data sequences---and typical data sequences have the feature that no computable forecasting system is high-low calibrated for them. Another purpose is to return with these results in mind to a question raised by \citet[\S 6]{Dawid:1982rc}---Is it really rational for Bayesian forecasters to be certain that they are calibrated?

\section{Meagre and Residual Sets and the Banach--Mazur Game}

In a probability space, the notions of sets of measure zero and sets of measure one give us natural ways of making precise the notions of sets so small as to be ignorable and of overwhelmingly large sets. In a topological space, the corresponding notions are as follows: a subset $Y$ of a topological space $X$ is \emph{meagre} (or \emph{first category}) if can be written as $Y=\bigcup_{k=1}^\infty Y_k,$ where each $Y_k$ is nowhere dense in $X$ (i.e., the closure of each $Y_k$ in $X$ has empty interior); the complement of a meagre subset of a topological space $X$ is called \emph{residual.} 

Taking as our prototype the notion of a countable subset of the real line, it is natural to think that any class of `small' subsets of a given space should be closed under taking subsets and countable unions. And it is natural to count a nowhere dense subset $Y$ of a topological space $X$ as small. This motivates the practice of topologists of treating properties as typical (rare) if they correspond to residual (meagre) subsets of a well-behaved space of interest. Thus, it is often said that typical continuous functions are nowhere differentiable because the set of such functions forms a residual subset of the space of continuous functions. 

The \emph{Banach--Mazur} game provides a convenient alternative characterization of meagre sets. Let $X$ be a topological space and let $\mathcal{G}$ be an open basis for the topology of  $X.$ Let a subset $A$ of $X$ be selected. A game is to be played. Player 1 and Player 2 will alternate taking turns picking elements $G_1,$ $G_2,$ $G_3,$ \ldots of $\mathcal{G}.$ The first move of Player 1 is unconstrained; but each subsequent $G_{k+1}$ must be a subset  of $G_k.$ Player 1 wins if $A$ has nonempty intersection with $\bigcap_{k=1}^\infty G_k.$ Otherwise Player 2 wins.  \citet{Oxtoby:1957ek} shows that $A$ is meagre in $X$ if and only if Player 2 has a winning strategy---no matter what Player 1 does, Player 2 can always win.

\section{Noncalibration is Typical}

For our forecasting systems, the space of possible data sequences is just the Cantor space $\mathcal{C}$---the space of infinite binary sequences. As a product space, this space has a natural topology: the topology of pointwise convergence. An open basis for this topology is given by the sets of the following form: for any finite binary string $w,$ $B_w$ is the set of binary sequences that have $w$ as an initial segment.

\begin{prop} Let $M$ be a forecasting system. Then relative to typical data sequences in $\mathcal{C},$ $M$ fails to be high-low calibrated.
\proof Let $A$ be the set of data sequences for which $M$ is high-low calibrated. We aim to show that $A$ is meagre in $\mathcal{C}.$ To this end,  we consider the Banach--Mazur game for $A\subset \mathcal{C},$ relative to the open basis for the topology of $\mathcal{C}$ consisting of sets of the form $B_w$ ($w$ a finite binary string). So Player 1 and Player 2 take turns playing binary strings, with with Player 1 winning if the sequence that results from concatenating the strings played lies in $A$, Player 2 winning otherwise. 
 
We need to show that Player 2 has a winning strategy. Here is how it works. Every time it is Player 2's turn, she plays a binary string, each bit of which is chosen to frustrate $M$---if Player 2 is playing the $k$th bit overall, then she plays a 1 if $\pi_k<.5,$ otherwise she plays a 0 for that bit. No matter what Player 1 has already done, by doing this for long enough, Player 2 can ensure that either: (i) the mean discrepancy for days on which the forecast probability is less than .5 is at least .25;  or (ii) the mean discrepancy for days on which the forecast probability is at least .5 is less than $-.25.$ Once she achieves one or the other of these objectives, she ends her turn. No matter how Player 1 plays, the infinite sequence $\boldsymbol{a}$ that results has the feature that it is fails to be calibrated either for the subsequence of days on which the forecast probability is less than .5 or for the subsequence of days on which the forecast probability is at least .5 (since at least (i) or (ii) must be achieved infinitely often by Player 2). Thus Player 2 has a winning strategy. 
\end{prop}
 So $M$ fails to be high-low calibrated for typical data sequences. Further, since the intersection of countably many residual sets is residual, and since there are only countably many computable forecasting systems, we see that typical data sequences have the feature that each computable forecasting system fails to be high-low calibrated for each of them.  
 
\section{Rationality, Calibration, and Arrogance}

Agents obeying the dictates of Bayesianism tend to be supremely self-confident. 

Consider the problem of reconstructing the i.i.d. probabilities governing a chance mechanism from knowledge of longer and longer finite sequences of its outputs. A Bayesian will posit a prior $\mu$ on the space $\Theta$ of parameter values under consideration, then update by conditionalization. Each prior assigns measure one to the subset of $\Theta$ consisting of those $\theta$ at which it is consistent (i.e., those $\theta$ such that the posterior probabilities resulting from conditionaliztion converge to a point mass at $\theta,$ with $\theta$-probability one). But these sets need not be small in any intuitive sense. Indeed, \citet{Freedman:1963it,Freedman:1965ph}  showed that when the chance process under consideration has a countable infinity of possible outcomes, typical priors fail to be consistent at typical $\theta$---although tail-free priors are consistent at every $\theta.$ It is often thought that there is a moral here: Bayesians should not settle for subjective certainty of success---priors inconsistent for vast swathes of possibilities should be viewed with suspicion.

The situation is a little worse with our forecasting problem: we again find that Bayesian agents are subjectively certain of success (i.e., of calibration)---but now we find that each and every one of them is typically unsuccessful. If we are to view poorly behaved priors with suspicion, then \emph{every} prior is under a cloud.
 
There is no special problem here for a pragmatic Bayesian statistician---the notion of calibration is remote from practical questions \citep{Schervish:1985ud}. But in some corners, at least, enthusiasm persists for the idea that rationality dictates that agents should have coherent subjective probabilities that they update by conditionalization. On this understanding, rationality \emph{requires} that we should be certain that we will be successful at certain tasks, even when we know that, among the possibilities under consideration, those at which we succeed are rare and those at which we fail are typical.

\bibliographystyle{plainnat}
	\bibliography{calibration}  
 \end{document}